\documentclass[11pt]{article}
\usepackage{amsmath}
\usepackage{dsfont}
\usepackage{mathrsfs}
\usepackage{amsmath,amssymb}
\usepackage{amsfonts}
\usepackage{hyperref}
\usepackage{amsthm}
\usepackage{graphicx}
\usepackage{subfigure}
\usepackage{xcolor}
\usepackage{cite}
\usepackage{epic,eepic,epsf,epsfig}
\usepackage{floatflt}
\hfuzz=\maxdimen
\theoremstyle{definition}
\newtheorem{lemma}{Lemma}[section]
\newtheorem{definition}[lemma]{Definition}
\newtheorem{proposition}[lemma]{Proposition}
\newtheorem{theorem}[lemma]{Theorem}

\newtheorem{remark}{Remark}

\numberwithin{equation}{section}

\DeclareFixedFont{\Acknowledgment}{OT1}{cmr}{bx}{n}{14pt}
\allowdisplaybreaks
\begin{document}

\title{\bf Combinatorial $p$-th Calabi flows on surfaces}
\author{Aijin Lin, Xiaoxiao Zhang}
\date{}
\maketitle

\begin{abstract}
For triangulated surfaces and any $p>1$, we introduce the combinatorial $p$-th Calabi flow which precisely equals the combinatorial Calabi flows first introduced in H. Ge's thesis \cite{Ge-thesis} (or see H. Ge \cite{Ge4}) when $p=2$. The difficulties for the generalizations come from the nonlinearity of the $p$-th flow equation when $p\neq 2$. Adopting different approaches, we show that the solution
to the combinatorial $p$-th Calabi flow exists for all time and converges if and only if there
exists a circle packing metric of constant (zero resp.) curvature in Euclidean (hyperbolic resp.) background geometry. Our results generalize the work of H. Ge \cite{Ge4}, Ge-Xu \cite{GX} and Ge-Hua \cite{GB} on the combinatorial Calabi flow from $p=2$ to any $p>1$.
\end{abstract}
\vspace{12pt}

\section{Introduction}
\hspace{14pt} Geometric flows are powerful tools to find canonical metrics on a given manifold which is an extremely important problem in differential geometry.
R. Hamilton \cite{RH} introduced the Ricci flow, which has been used to prove the uniformization theorem and solve the Poincar$\acute{\textrm{e}}$ conjecture. Seeking constant curvature metrics, E. Calabi \cite{Ca, Ca2} introduced the Calabi flow and studied the variational
problem of minimizing the so-called ¡°Calabi energy¡± in any fixed cohomology class
of K$\ddot{\text{a}}$hler metrics. For dimension
two, it is proved that both the Calabi flow and the
normalized Ricci flow exist for all time and converge to a constant scalar curvature
metric (cf. \cite{Cha, Cha2, Chen} for further references).\par
Suppose $X$ is a closed surface with a triangulation
$\mathcal{T} = (V,E,F)$, where $V , E, F$ denote the sets of vertices, edges, and
faces respectively. Throughout this paper, a function defined on vertices is
an $N$-dimensional column vector, where $N = |V |$ is the number of vertices. Moreover,
all vertices, marked by $v_{1}, . . . , v_{N}$, are ordered one by one.
\begin{definition}\label{CP}
A \emph{circle packing metric} is defined to be a positive function
 \begin{align*}
 r : V&\rightarrow(0,+\infty)\\
     v_{i}&\mapsto r_{i}, i = 1, . . . , N
 \end{align*}
  on the vertices.
\end{definition}
 Hence we may think of circle packing metrics
as points in $\mathbb{R}^{N}_{>0}$, where $\mathbb{R}^{N}_{>0}$ means $N$ times of Cartesian product of $(0,+\infty)$. A \emph{weight} on the triangulation is defined to be a
function $\Phi : E\rightarrow[0, ¦Ð/2]$.  A
triangulated surface with a weight $\Phi$ is denoted as $(X, T, \Phi)$.
Let $l : E\rightarrow(0,+\infty)$ be a positive function assigning each edge ${i, j} \in E$ a
length $l_{ij}$. $l$ is called a \emph{piecewise linear metric} if for every triangle ${i, j, k} \in F$,
the three edge lengthes $l_{ij}, l_{jk}, l_{ik}$ satisfy triangle inequalities. For a fixed
triangulated surface $(X, \mathcal{T}, \Phi)$, every circle packing metric $r$ determines a piecewise
linear metric on $X$ by setting the length of edge ${i, j} \in E$ as
$$l_{ij} =\sqrt{r_{i}^{2}+r_{j}^{2}+2r_{i}r_{j}\cos(\Phi_{ij})},$$
in Euclidean background geometry and
$$l_{ij} =\cosh^{-1}(\cosh r_{i}\cosh r_{j}+\sinh r_{i}\sinh r_{j}\cos(\Phi_{ij})),$$ in hyperbolic background geometry.
As a consequence, each face in $F$ is isometric to a Euclidean triangle (hyperbolic resp.). Specifically,
each face ${i, j, k} \in F$ is a Euclidean triangle (hyperbolic resp.) with edge lengths $l_{ij}, l_{jk}, l_{ik}$
because $l_{ij}, l_{jk}, l_{ik}$ satisfy triangle inequalities (Lemma 13.7.2 \cite{T}). Furthermore,
the triangulated surface $(X, \mathcal{T})$ is composed by gluing Euclidean triangles coherently (hyperbolic resp.).\par
 Given a
triangulated surface $(X, \mathcal{T}, \Phi)$ with a circle packing metric $r$, all inner angles of the
triangles are determined by $r_{1}, . . . , r_{N}$. Denote $\theta_{i}^{jk}$
as the inner angle at vertex $i$ in
the triangle ${i, j, k} \in F$. The well-known \emph{combinatorial Gauss curvature} $K_{i}$ at vertex $i$ is defined as
\begin{equation}
K_{i} = 2\pi-\sum_{{i, j, k} \in F}\theta_{i}^{jk},
\end{equation}
where the sum is taken over each triangle with $i$ as one of its vertices.
For every circle packing metric $r$ on $(X, \mathcal{T}, \Phi)$, we have the following \emph{combinatorial Gauss-Bonnet formula} introduced by Thurston in \cite{T}.
\begin{equation}\label{GB}
\sum_{i=1}^{N}K_{i}=2\pi\chi(X)-\lambda \text{Area}(X),
\end{equation}
where $\lambda=-1, 0, 1$ correspond three geometries, i.e., hyperbolic geometry $H^{2}$, Euclidean geometry $E^{2}$ and spherical geometry $S^{2}$.\par

\begin{definition}\label{Cons}
A \emph{constant curvature circle packing metric} $r_{av}$ is a circle packing metric that
determines the constant combinatorial curvature $K_{av}= K(r_{av}) = k_{av}(1, . . . , 1)^{T}$ , where $k_{av} = 2¦Ð\chi(X)/N$.
\end{definition}
The constant (zero resp.) curvature circle packing metric is
a good candidate for privileged metrics. Thurston \cite{T} first studied this class of metrics,
and found that there are combinatorial obstructions for the existence of constant (zero resp.)
curvature metrics . Given a triangulated surface $(X, \mathcal{T}, \Phi)$, Thurston introduced the circle packing metric and found that there are combinatorial obstructions for the existence of a circle
packing metric with constant (zero resp.) combinatorial curvatures (see section 13.7 in \cite{T}). \par
Later, Chow-Luo \cite{CL} introduced the
\emph{combinatorial Ricci flow in Euclidean background geometry}
\begin{equation}\label{R1-flow}
  \frac{dr_{i}}{dt}=-K_{i}r_{i},
\end{equation}
and the \emph{normalized combinatorial Ricci flow}
\begin{equation}\label{R2-flow}
  \frac{dr_{i}}{dt}=(k_{av}-K_{i})r_{i},
\end{equation}
 which is equivalent to (\ref{R1-flow}). The
\emph{combinatorial Ricci flow in hyperbolic background geometry} is defined by Chow-Luo \cite{CL} as
\begin{equation}\label{Rh-flow}
  \frac{dr_{i}}{dt}=-K_{i}\sinh r_{i}.
\end{equation}
Then Chow-Luo proved that the combinatorial Ricci flows above exist for
all time and converge exponentially fast to Thurston's circle packing on surfaces if and only if there
exists a circle packing metric of constant (zero resp.) curvature in Euclidean (hyperbolic resp.) background geometry.
Since then, various discrete curvature flows were introduced and studied. We refer the the readers to H. Ge \cite{Ge1, Ge2, Ge3}, Ge-Hua-Jiang \cite{GBJ} Ge-Jiang \cite{GJ1, GJ2, GJ3}, Zhang-Chang \cite{ZC}, Lin-Zhang \cite{LZ} and Zhang-Lin \cite{ZL, ZL2}.
Specially, motivated by the work of Chow-Luo, H. Ge in his thesis \cite{Ge-thesis} (or see \cite{Ge4}\cite{GB}\cite{GX}) introduced the \emph{combinatorial Calabi flow in Euclidean background geometry}
\begin{equation}\label{C1-flow}
 \frac{dr_{i}}{dt}=\Delta K_{i}r_{i}=-(L^{T}K)_{i}r_{i},
\end{equation}
and the \emph{combinatorial Calabi flow in hyperbolic background geometry}
\begin{equation}\label{C2-flow}
 \frac{dr_{i}}{dt}=\Delta K_{i}\sinh r_{i}=-(L^{T}K)_{i}\sinh r_{i},
\end{equation}
where
$\Delta K_{i}=-\sum_{j\sim i}L_{ji}K_{j}$, $L_{ij}=\frac{\partial K_{i}}{\partial r_{j}}r_{j}$ in Euclidean background geometry, and $L_{ij}=\frac{\partial K_{i}}{\partial r_{j}}\sinh r_{j}$ in hyperbolic background geometry. \par
 Similar to the combinatorial Ricci flows, H. Ge \cite{Ge4}, Ge-Xu \cite{GX} and Ge-Hua \cite{GB} proved the combinatorial Calabi flows have the long time existence and converge exponentially fast to Thurston's circle packing on surfaces if and only if there
exists a circle packing metric of constant (zero resp.) curvature in Euclidean (hyperbolic resp.) background geometry.

 In this paper, we generalize above results on the combinatorial Calabi flow for $p = 2$ to any $p>1$.  The difficulties for the generalizations mainly come from the fact that the $p$-th flow equations are nonlinear when $p\neq 2$. First for any $p>1$ we introduce the combinatorial $p$-th Calabi flow.
\begin{definition}
The \emph{combinatorial $p$-th Calabi flow in Euclidean background geometry} is defined as
\begin{equation}\label{p-Ca}
\frac{dr_{i}}{dt}=\Delta_pK_{i}r_{i},
\end{equation}
where $\Delta_p$ is the \emph{discrete Laplace operator} defined as
\begin{equation}\label{p-Lap}
  \Delta_pf_i=\sum_{j\thicksim i}B_{ij}|f_j-f_i|^{p-2}(f_j-f_i), \forall f\in \mathbb{R}^{N}
\end{equation}
with $$B_{ij}=\frac{\partial(\theta_{i}^{jk}+\theta_{i}^{jl})}{\partial r_{j}}r_{j}.$$
Here $\theta_{i}^{jk}$ is the inner angle at the vertex $v_{i}$ in a triangle $\triangle ijk\in F$.\par
The \emph{combinatorial $p$-th Calabi flow in hyperbolic background geometry} is defined as
\begin{equation}\label{p-Cah}
\frac{dr_{i}}{dt}=(\Delta_pK_{i}-A_{i}K_{i})\sinh r_{i},
\end{equation}
where $A_{i}$ is defined as
\begin{equation*}
 A_{_i}=\sinh r_i\frac{\partial}{\partial r_i}\Big(\sum_{\{ijk\}\in F} \text{Area}(\triangle v_iv_jv_k )\Big)
\end{equation*}
\end{definition}
\begin{remark}
When $p=2$, by Proposition \ref{L} and (\ref{LB}) (see Section 2 for details), the combinatorial $p$-th Calabi flows (\ref{p-Ca}) and (\ref{p-Cah}) are exactly the combinatorial Calabi flows (\ref{C1-flow}) and (\ref{C2-flow}) introduced by H. Ge \cite{Ge-thesis}.\par
\end{remark}
Now we state our main results in this paper as follows.\par

\begin{theorem}\label{Ca-E}
\emph{Suppose $(T, \Phi)$ is a weighted generalized triangulation of a closed connected surface $X$. Given any initial circle packing
metric $r(0) \in \mathbb{R}^{N}_{>0}$ satisfying $\prod_{i=1}^{N}r_{i}(0)=1$, the solution to the combinatorial $p$-th Calabi flow \ref{p-Ca} in Euclidean background geometry exists for $t\in
[0, +\infty)$. Additionally, $r(t)$ converges if and only if there exists a constant curvature circle packing metric $r_{av}$}.
\end{theorem}

\begin{theorem}\label{Ca-H}
\emph{Suppose $(T, \Phi)$ is a weighted generalized triangulation of a closed connected surface $X$ of negative Euler characteristic. Given any initial circle packing
metric $r(0) \in \mathbb{R}^{N}_{>0}$, the solution to the combinatorial $p$-th Calabi flow \ref{p-Cah} in hyperbolic background geometry exists for $t\in
[0, +\infty)$. Additionally, $r(t)$ converges if and only if there exists a zero curvature circle packing metric $r_{ze}$}.
\end{theorem}
Theorem \ref{Ca-E} generalizes H. Ge's results \cite{Ge4} on the combinatorial Calabi flow in Euclidean background geometry.
 Theorem \ref{Ca-H} generalizes the results of Ge-Xu \cite{GX} and Ge-Hua \cite{GB} on the combinatorial Calabi flow in hyperbolic background geometry. Especially, different from the the combinatorial Calabi flow ($p=2$), our results show that the combinatorial $p$-th Calabi flow does not have the exponential convergence for $p>1, p\neq2$. \par

The key point of the proof of Theorem \ref{Ca-E} is that we succeed showing the compactness of the solution $\{r(t)\}$ to the combinatorial $p$-th Calabi flow (\ref{p-Ca}). More precisely, in the case of Euclidean background geometry, similar to the combinatorial Ricci flow and the combinatorial Calabi flow, we show that the solution $\{r(t)\}$ to the combinatorial $p$-th Calabi flow is always in a hypersurface. When restricted to the hypersurface, the combinatorial Ricci potential (see Section 2 for the definition) first introduced by de Verdiere \cite{DV} satisfies the Palais-Smale condition and further is proper, which suffices to establish the the compactness of the solution $\{r(t)\}$ to the combinatorial $p$-th Calabi flow. In the case of hyperbolic background geometry, we also show the long time existence and the uniform convergence. The difficulty of the proof of Theorem \ref{Ca-H} comes from the long time existence in this case. By generalizing the work of Ge-Xu \cite{GX} and Ge-Hua \cite{GB}, we can show that the boundness of the solution $\{r(t)\}$ to the combinatorial $p$-th Calabi flow after getting rid of the additional assumption in \cite{GX} on curvature.\par

This paper is organized as follows. In Section 2, we give some preliminaries.
In Section 3, we establish a uniform estimate which plays an important role in showing the long time existence of the solution to the combinatorial $p$-th Calabi flow in hyperbolic background geometry. Then in Section 4, we prove the long time existence parts of both Theorem \ref{Ca-E} and Theorem \ref{Ca-H}. In Section 5, we prove the long time convergence parts of both Theorem \ref{Ca-E} and Theorem \ref{Ca-H}. In section 6, we consider other combinatorial $p$-th curvature flows and discuss some interesting unsolved questions for future study.\\

\noindent{\bf Acknowledgments.}
The first author would like to thank Professor Gang Tian and Professor Huijun Fan for constant encouragement and support. The second author would like to thank Professor Yanxun Chang for constant guidance and encouragement. Both authors would also like to thank Professor Huabin Ge for many helpful conversations. The first author is supported by National Natural Science Foundation of China under Grant No. 11401578. The second author is supported by National Natural Science Foundation of China under Grant No. 11471138, No. 11501027 and No. 11871094 and Fundamental Research Funds for the Central Universities Grant No. 2017JBM072.
\section{Preliminaries}
\hspace{14pt}  In this section, we give necessary preliminaries.
\subsection{Some basic propositions}
\hspace{14pt} Any circle packing metric
$r$ determines an intrinsic metric structure on fixed $(X, T, \Phi)$ by Euclidean cosine
law. The lengths $l_{ij}$ , angles $\theta^{jk}_{i}$
 and curvatures $K_{i}$ are elementary functions of
$r = (r_{1}, . . . , r_{N})^{T}$. \par
Let the background geometry $K^{2}$ be one of the three geometries, $E^{2}$ ( Euclidean ), $H^{2}$ ( hyperbolic ), and $S^{2}$ ( spherical ), i.e., a simply connected surface with a complete Riemannian metric of curvature equal to 0, 1, and -1. For convenience, we choose coordinate transformations according to the background geometry. To be more precise, we choose $u_{i}=\ln r_{i}$ when $K^{2}=E^{2}$, $u_{i}=\text{ln }\tanh (r_{i}/2)$ when $K^{2}=H^{2}$ and $u_{i}=\ln \tan (r_{i}/2)$ when $K^{2}=S^{2}$.\par
Denote $j\sim i$ if the vertices $i$ and $j$ are adjacent. For any
vertex $i$ and any edge $j\sim i$, in the case of $E^{2}$, set
\begin{equation}\label{B1}
  B_{ij}=\frac{\partial(\theta_{i}^{jk}+\theta_{i}^{jl})}{\partial r_{j}}r_{j}..
\end{equation}

In the case of $H^{2}$, set
\begin{equation}\label{B2}
 B_{ij}=\frac{\partial(\theta_{i}^{jk}+\theta_{i}^{jl})}{\partial r_{j}}\sinh r_{j},
\end{equation}
\begin{equation}\label{A}
 A_i=\sinh r_i\frac{\partial}{\partial r_i}\Big(\sum_{\{ijk\}\in F} \text{Area}(\triangle v_iv_jv_k )\Big),
\end{equation}
where $k,l$ are the vertices such that $\{ijk\},\ \{ijl\}$ are adjacent faces. \par
Note $B_{ij} = B_{ji}$ since $\frac{\partial \theta_{i}^{jk}}{\partial r_{j}}r_{j}=\frac{\partial \theta_{j}^{ik}}{\partial r_{i}}r_{i}$ (see Lemma 2.3 in \cite{CL}).
We have the following estimate for $B_{ij}$ (see Lemma 2.2 in \cite{CL} or Proposition 3.1 in \cite{Ge4})

\begin{proposition}\label{B}
\emph{When $K^{2}=H^{2}$, then for any $1\leq i,j\leq N$ and $i\sim j$, we have
\begin{equation}
  B_{ij}>0.
\end{equation}
When $K^{2}=E^{2}$, then for any $1\leq i,j\leq N$ and $i\sim j$, we have
\begin{equation}
  0<B_{ij}<2\sqrt{3}.
\end{equation}}
\end{proposition}

Set $A=diag\{A_{1}, \cdots, A_{N}\}$ and $ L_{B}=((L_{B})_{ij})_{N \times N}$, where
\begin{eqnarray}\label{LB}
 (L_{B})_{ij}=\left\{
\begin{aligned}
&\sum_{k\sim i}B_{ik}, & j=i, \\
&-B_{ij}, & j\sim i, \\
&0, &\text{else}.
\end{aligned}
\right.
\end{eqnarray}
Recall $L=(L_{ij})=(\frac{\partial K_{i}}{\partial u_{j}})$, then we have the following proposition (see Theorem 3.1 in \cite{GX}).
\begin{proposition}\label{L}
\emph{In the case of $E^{2}$, we have $L=L_{B}$. In the case of $H^{2}$, we have $L=A+L_{B}$. Furthermore,
$L_{B}$ is semi-positive definite, and $A$ is positive definite}.
\end{proposition}

As a corollary, we have the following proposition (see Proposition 3.3 in \cite{Ge4}).

\begin{proposition}\label{L2}
\emph{If $K^{2}=H^{2}$, then $L$ is positive definite.
If $K^{2}=E^{2}$, then $L$ is a semi-positive definite $N \times N$ matrix, whose rank is $N-1$.
Moreover, the null space of $L$ is $\text{Ker}(L) = \{t(1, . . . , 1)^{T} |t \in \mathbb{R}\}$}.
\end{proposition}

\subsection{Combinatorial Ricci potential}
For $K^{2}=E^{2}$ or $H^{2}$, we have $\frac{\partial K_{i}}{\partial u_{j}}=\frac{\partial K_{j}}{\partial u_{i}}$ by $\frac{\partial \theta_{i}^{jk}}{\partial r_{j}}r_{j}=\frac{\partial \theta_{j}^{ik}}{\partial r_{i}}r_{i}$ (see Lemma 2.3 in \cite{CL}). Thus the smooth differential 1-form $\sum_{i=1}^{N}(K_i-k_{av})du_i$ is closed, so we can consider the combinatorial Ricci potential.
\begin{definition}\label{R-pot}
When $K^{2}=E^{2}$, the following function
\begin{equation}\label{RP2}
F(u)\triangleq\int_{a}^u\sum_{i=1}^{N}(K_i-k_{av})du_i,\ u\in\mathbb{R}^{N}.
\end{equation}
is called \emph{combinatorial Ricci potential}, where $a$ is any point in $\mathbb{R}^{N}$. In the case of $H^{2}$ the \emph{combinatorial Ricci potential} is defined as
\begin{equation}\label{RP1}
F(u)\triangleq\int_{a}^u\sum_{i=1}^{N}K_idu_i,\ u\in\mathbb{R}_{<0}^{N},
\end{equation}
where $a$ is any point in $\mathbb{R}_{<0}^{N}$.
\end{definition}
These integrals were first introduced by de Verdiere \cite{DV}, who also first proved the convexity of $F$ in the zero weight case. Chow-Luo proved the following convex theorem (see Proposition 3.9 in \cite{CL}).
\begin{proposition}\label{convex}
\emph{When $K^{2}=H^{2}$, then the combinatorial Ricci potential (\ref{RP1}) $F(u): \mathbb{R}_{>0}^{N}\rightarrow \mathbb{R}$ is strictly convex. When $K^{2}=E^{2}$, then (\ref{RP2}) $F(u): \mathbb{R}^{N}\rightarrow \mathbb{R}$ is convex. Furthermore, it satisfies $F(u+c(1, 1, \cdots, 1)^{T})=F(u)$ for any $a\in\mathbb{R}^{N}$ and $F$ is strictly convex when restricted in the hyperplane $\mathcal{U}_{0}=\{(u_{1}, \cdots, u_{N})^{T}\mid \sum_{i=1}^{N}u_{i}=0\}$}.
\end{proposition}

As a consequence of the above proposition, we have the following rigidity theorem (see Corollary 3.11 in \cite{CL}).
\begin{proposition}\label{rigidity}
\emph{Let $\Pi: \mathbb{R}_{>0}^{N}\rightarrow \mathbb{R}^{N}$ be the map sending a metric $r_{1}, \cdots, r_{N})^{T}$ to the corresponding curvature $(K_{1}, \cdots, K_{N})^{T}$ where $K^{2}=H^{2}$ or $E^{2}$. Then \\
($\text{a}$) In the case of $H^{2}$, $\Pi$ is injective, i.e., the metric is determined by its curvature.\\
($\text{b}$) In the case of $E^{2}$, $\Pi$ restricted to the hypersurface $\mathcal{P}_{1}=\{(r_{1}, \cdots, r_{N})^{T}\in  \mathbb{R}_{>0}^{N}\mid \prod_{i=1}^{N}r_{i}=1\}$ is injective, i.e., the metric is determined by its curvature up to a scalar multiplication}.
\end{proposition}
To prove our main results, it is key to show that $F$ satisfies the Palais-Smale condition. In fact, Chow-Luo \cite{CL} has proved the following theorem
\begin{proposition}\label{PS}
\emph{If $K^{2}=H^{2}$, then the combinatorial Ricci potential $F$ (\ref{RP1}) satisfies the Palais-Smale condition if the following Thurston's combinatorial conditions hold:\\
$(\textbf{c1})$ for any three edges $e_1,e_2,e_3$ forming a null homotopic loop in $M$, if
$\sum_{i=1}^{3}\Phi(e_i)\geq\pi$, then $e_1,e_2,e_3$ form the boundary of a triangle in $F$;\\
$(\textbf{c2})$ for any four edges $e_1,e_2,e_3,e_4$ forming a null homotopic loop in $M$, if $\sum_{i=1}^{4}\Phi(e_i)\geq2\pi$, then $e_1,e_2,e_3,e_4$ form the boundary of the union of two adjacent triangles.\\
If $K^{2}=E^{2}$, then the combinatorial Ricci potential $F$ (\ref{RP2}) restricted to the hyperplane $\mathcal{U}_{0}=\{(u_{1}, \cdots, u_{N})^{T}\mid \sum_{i=1}^{N}u_{i}=0\}$, i.e., $F|_{\mathcal{U}_{0}}$ satisfies the Palais-Smale condition if for any proper subset $I\subset V$},
\begin{equation}\label{E-con}
2\pi|I|\chi(X)/N>-\sum_{(e,v)\in \text{Lk}(I)}(\pi-\Phi(e))+2\pi\chi(F_{I}),
\end{equation}
\emph{where $F_{I}$ is the set of all cells in $T$ whose vertices are in $I$, and $\text{Lk}(I)$ is the link of $I$ which is the set of pairs $(e, v)$ of an edge $e$ and a vertex $v$ so that (1) the end points of $e$ are not in $I$ and (2) the vertex $v$ is in $I$ and (3) $e, v$ form a triangle}.
\end{proposition}
\begin{remark}
By the work of Thurston \cite{T} on circle packing, conditions (c1) and (c2) are equivalent to the existence of the zero curvature circle packing metric when $K^{2}=H^{2}$, and condition (\ref{E-con}) is equivalent to the existence of the constant curvature circle packing metric $K^{2}=E^{2}$.
\end{remark}

Furthermore, we can show the properness of $F$ under conditions above.  In fact, H. Ge and Ge-Xu have proved the following theorem (see Theorem B.2 in \cite{Ge4} and Lemma B.1 in \cite{GX}).
\begin{proposition}\label{proper}
\emph{Making the same assumptions as in Proposition \ref{PS}, then for the combinatorial Ricci potential $F$ (\ref{RP2}) and (\ref{RP1})
\begin{equation}\label{prop}
\lim_{\|u\|\rightarrow+\infty}F(u)=+\infty.
\end{equation}
Moreover, the combinatorial Ricci potential $F$ (\ref{RP1}) is proper in the case of $H^{2}$. The combinatorial Ricci potential $F$ (\ref{RP2}) is proper when restricted to the hyperplane $\mathcal{U}_{0}$ in the case of $E^{2}$}.
\end{proposition}

\subsection{Combinatorial $p$-th Calabi flows}
Recall the combinatorial $p$-th Calabi flow in the case of $E^{2}$ is
\begin{equation*}
\frac{dr_{i}}{dt}=\Delta_pK_{i}r_{i},
\end{equation*}
where
\begin{equation*}
  \Delta_pK_i=\sum_{j\thicksim i}B_{ij}|K_j-K_i|^{p-2}(K_j-K_i),
\end{equation*}
and the combinatorial $p$-th Calabi flow in the case of $H^{2}$ is
\begin{equation*}
\frac{dr_{i}}{dt}=(\Delta_pK_{i}-A_{i}K_{i})\sinh r_{i},
\end{equation*}
where
\begin{equation*}
 A_{_i}=\sinh r_i\frac{\partial}{\partial r_i}\Big(\sum_{\{ijk\}\in F} \text{Area}(\triangle v_iv_jv_k )\Big)
\end{equation*}
By coordinate transformations, more precisely, $u_{i}=\ln r_{i}$ when $K^{2}=E^{2}$, and $u_{i}=\ln \tanh (r_{i}/2)$ when $K^{2}=H^{2}$,  then we can rewrite the combinatorial $p$-th Calabi flow as
\begin{equation}\label{p-Ca-u}
 u'_{i}(t)=\Delta_pK_{i},
\end{equation}
in the case of $E^{2}$ and the combinatorial $p$-th Calabi flow
\begin{equation}\label{p-Cah-u}
u_{i}'(t)=L_pK_{i}=\Delta_pK_{i}-A_{i}K_{i},
\end{equation}
in the case of $H^{2}$.\par
We can also write the combinatorial $p$-th Calabi flow equations in a matrix
form as
\begin{equation}\label{p-flow-m-e}
  u'(t)=\Delta_pK
\end{equation}
in the case of $E^{2}$  and
\begin{equation}\label{p-flow-m-h}
  u'(t)=L_pK=(\Delta_p-A)K
\end{equation}
in the case of $H^{2}$,
where $A=diag\{A_1,A_2,\ldots,A_N\}$. The matrix $A$ is positive definite by Proposition \ref{L}.\par

\section{A uniform estimate}

\hspace{14pt}  In order to prove Theorem \ref{Ca-H}, we need the following uniform estimate which plays an important role in our arguments. First recall a lemma (see Lemma 13.7.3 in \cite{T} ) proved by Thurston  as follows:
\begin{lemma}\label{Thurston}
\emph{ Let $\triangle v_iv_jv_k$ be a hyperbolic triangle which is patterned by three circles with fixed weights $\Phi_{ij},\Phi_{jk},\Phi_{ik}\in[0,\pi/2]$ as intersection angles. Let $\theta_i^{jk}$ be the inner angle at $v_i$. Then $\partial \theta_i^{jk}/\partial r_j>0,\ \partial \theta_i^{jk}/\partial r_i<0$ and $\partial (\theta_i^{jk}+\theta_j^{ik}+\theta_k^{ij})/\partial r_i<0$}.
\end{lemma}
Generalizing Lemma 3.2 in \cite{GB}, we have the following uniform estimate which plays an important role in our arguments.
\begin{lemma}\label{estimate}
\emph{For any $\lambda>1$, there exists a constant $C$, depending on $\lambda$, such that if $r_i\geq C$, then}
\begin{equation}\label{hh}
  \frac{\partial }{\partial r_i}(\lambda \text{Area}(\triangle v_iv_jv_k)+\theta_i^{jk})\geq 0.
\end{equation}
\end{lemma}

\noindent\begin{proof}The proof is similar to the proof of Lemma 3.2 in \cite{GB}.
Adopt the argument from the hyperbolic geometry. Suppose that the triangle
 $\triangle v_iv_jv_k$ is embedded in the half hyperbolic plane $\mathbb{H}^2$, with $v_j, v_k$ and the corresponding radii $r_j, r_k$ fixed. Let $\hat{v}_i$
be the new vertex with a larger radius $\hat{r}_i>r_i$. Then $l_{\hat{v}_iv_j} > l_{v_iv_j}$ and
$l_{\hat{v}_iv_k} > l_{v_iv_k}$ due to $\hat{r}_i>r_i$. Draw two triangles, $\triangle v_iv_jv_k$ and $\triangle \hat{v}_iv_jv_k$ whose common edge $v_jv_k$ is in the
same half hyperbolic plane separated by the (extended) geodesic $v_jv_k$. By Lemma \ref{Thurston}, for
fixed $r_j$ and $r_k$ the angles $\theta_j^{ik}$ and $\theta_i^{jk}$
are increasing which implies that the vertex $v_i$
lies in the interior of the triangle $\triangle v_iv_jv_k$ (please refer to Figure 2. in \cite{GB}). Denote by $\hat{\theta}_i^{jk}$, $\hat{\theta}_j^{ik}$ and $\hat{\theta}_k^{ij}$
 three inner angles of the new triangle $\triangle v_iv_jv_k$ respectively. To prove (\ref{hh}), it is sufficient to show that for any $\hat{r}_i>r_i$ which is sufficiently close to $r_i$,
 \begin{equation}\label{equ-h}
   \lambda\text{Area}(\triangle \hat{v}_iv_jv_k)-\lambda\text{Area}(\triangle \hat{v}_iv_jv_k)+\hat{\theta}_i^{jk}-\theta_i^{jk}\geq0.
 \end{equation}
Set $x=\hat{\theta}_j^{ik}-\theta_j^{ik}$ and $\hat{\theta}_k^{ij}-\theta_k^{ij}$, then we have
\begin{align}
&\lambda\text{Area}(\triangle \hat{v}_iv_jv_k)-\lambda\text{Area}(\triangle \hat{v}_iv_jv_k)+\hat{\theta}_i^{jk}-\theta_i^{jk} \nonumber\\=&(\lambda-1)\text{Area}(\triangle \hat{v}_iv_jv_k)-(\lambda-1)\text{Area}(\triangle \hat{v}_iv_jv_k)-x-y\nonumber \\=&[(\lambda-1)\text{Area}(\triangle \hat{v}_iv_jv_k)-x]+[(\lambda-1)\text{Area}(\triangle \hat{v}_iv_jv_k)-y]
\end{align}
Then it suffices to prove the following inequalities
\begin{equation*}
  (\lambda-1)\text{Area}(\triangle \hat{v}_iv_jv_k)\geq x,\ \text{and}\ (\lambda-1)\text{Area}(\triangle \hat{v}_iv_jv_k)\geq y.
\end{equation*}
Without loss of generality, by the symmetry, we show that
\begin{equation}\label{suff-p}
  (\lambda-1)\text{Area}(\triangle \hat{v}_iv_jv_k)\geq x.
\end{equation}
Let $s$ be the point which attains the minimum distance from the
vertex $v_i$ to a point on the geodesic $\hat{v}_iv_j$, then $s$ is in the interior of the geodesic $\hat{v}_iv_j$ (please refer to Figure 2. in \cite{GB}) since $l_{\hat{v}_iv_j} > l_{v_iv_j}$. Suppose that $\hat{r}_i$ is sufficiently close to $r_i$ such that $l_{\hat{v}_iv_i}\leq1$. By the hyperbolic cosine law,
$$\cos x=\frac{\cosh l_{v_iv_j}\cosh l_{\hat{v}_iv_j}-\cosh l_{\hat{v}_iv_i}}{\sinh l_{v_iv_j}\sinh l_{\hat{v}_iv_j}}\rightarrow 1,$$
uniformly as $r_i\rightarrow \infty$. Hence there is a universal constant $C_1$ such that if $r_i\geq C_1$, then \begin{equation}\label{x-e}
x\leq\frac{\pi(\lambda-1)}{4\lambda}.
\end{equation}
Set $\beta=\angle sv_iv_j$, then there are two possibilities.\par
If $\beta<\frac{\pi}{4}$, then $\beta+x+\frac{x}{\lambda-1}<\frac{\pi}{2}$  by (\ref{x-e}). Considering the Gauss-Bonnet theorem in hyperbolic background geometry, i.e.,
$\beta+x+\frac{\pi}{2}=\pi-\text{Area}(\triangle sv_iv_j),$
we have $\text{Area}(\triangle \hat{v}_iv_iv_j)\geq\text{Area}(\triangle sv_iv_j)>\frac{x}{\lambda-1},$ which leads to (\ref{suff-p}).\par

If $\beta\geq\frac{\pi}{4}$, then
$$\frac{\sinh l_{sv_j}}{\sinh l_{v_iv_j}}=\sin \beta\geq\frac{\sqrt2}{2}.$$
For the hyperbolic right triangle $\triangle sv_iv_j$, using the hyperbolic cosine law, we have
$$\cos \beta=\sin x\cosh l_{sv_j},\ \cos x=\tanh l_{sv_j}/\tanh l_{v_iv_j},\ \sin \beta=\sinh l_{sv_j}/\sinh l_{v_iv_j}$$
This yields
\begin{align}\label{sin}
\sinh(\text{Area}(\triangle sv_iv_j))&=\cos(x+\beta)=\cos x\cos\beta-\sin x\sin\beta\nonumber\\&=\frac{\tanh l_{sv_j}}{\tanh l_{v_iv_j}}\sin x\cosh l_{sv_j}-\sin x\frac{\sinh l_{sv_j}}{\sinh l_{v_iv_j}}\nonumber\\&=\sin x(\cosh l_{v_iv_j}-1)\frac{\sinh l_{sv_j}}{\sinh l_{v_iv_j}}\nonumber\\&\geq\sin x\frac{\cosh r_i-1}{\sqrt 2}.
\end{align}

When $\lambda\geq2$, by (\ref{sin}) there exists a universal constant $C_2$ such that if $r_i\geq C_2$,
$$\sin(\text{Area}(\triangle sv_iv_j))\geq\sin x.$$
Noting that both $x$ and $\text{Area}(\triangle sv_iv_j)$ are in $(0,\frac{\pi}{2})$, we have $\text{Area}(\triangle sv_iv_j)\geq x$, which yields $(\lambda-1)\text{Area}(\triangle sv_iv_j)\geq\text{Area}(\triangle sv_iv_j)\geq x$.

When $1<\lambda<2$, by (\ref{sin}), we have
$$\sin((\lambda-1)\text{Area}(\triangle sv_iv_j))\geq(\lambda-1)\sin(\text{Area}(\triangle sv_iv_j))\geq\sin x\frac{\cosh r_i-1}{\sqrt 2}(\lambda-1).$$
Note that both $x$ and $(\lambda-1)\text{Area}(\triangle sv_iv_j)$ are in $(0,\frac{\pi}{2})$, thus there exists a a universal constant $C_3$, depending on $\lambda$, such that if $r_i\geq C_3$,
$$(\lambda-1)\text{Area}(\triangle sv_iv_j)\geq x.$$
$\hfill\Box$
\end{proof}

\begin{proposition}\label{Area-e}
\emph{For any $\epsilon>0$, there exists a universal number $C>0$, such that if $r_i>C$, then}
$$\epsilon A_i\geq\sum_{j\thicksim i}B_{ij}.$$
\end{proposition}
\noindent\begin{proof}
By (\ref{hh}), we have
\begin{align*}
&\epsilon A_i-\sum_{j\thicksim i}B_{ij}\\&=\epsilon\sum_{\{ijk\}\in F}\frac{\partial \text{Area}(\triangle v_iv_j v_k)}{\partial r_i}\sinh r_i-\sum_{j\thicksim i}\left(\frac{\partial\theta_i^{jk}}{\partial r_j}\sinh r_j+\frac{\partial\theta_i^{jl}}{\partial r_j}\sinh r_j\right)
\\&=\epsilon\sum_{\{ijk\}\in F}\frac{\partial \text{Area}(\triangle v_iv_j v_k)}{\partial r_i}\sinh r_i-\sum_{\{ijk\}\in F}\left(\frac{\partial\theta_i^{jk}}{\partial r_j}\sinh r_j+\frac{\partial\theta_i^{jl}}{\partial r_j}\sinh r_j\right)
\\&=\sum_{\{ijk\}\in F}\frac{\partial\left(\epsilon \text{Area}(\triangle v_iv_j v_k)-\theta_j^{ik}-\theta_k^{ij}\right)}{\partial r_i}\sinh r_i\\&=\sum_{\{ijk\}\in F}\frac{\partial\left((\epsilon+1) \text{Area}(\triangle v_iv_j v_k)+\theta_i^{jk}\right)}{\partial r_i}\sinh r_i\geq0.\quad\quad\Box
\end{align*}
\end{proof}

\section{The long time existence}
\hspace{14pt}  In this section we prove the long time existence parts of both Theorem \ref{Ca-E} and Theorem \ref{Ca-H}. In fact, we get the following theorem.
\begin{theorem}\label{p-ca-exis}
\emph{For any initial circle packing metric $r(0)\in\mathbb{R}_{>0}^{N}$, both the solutions to the combinatorial $p$-th Calabi flows (\ref{p-Ca}) in $E^{2}$ and (\ref{p-Cah}) in $H^{2}$ exist for all time $t\in [0, +\infty)$}.
\end{theorem}
Before showing the long time existence of combinatorial $p$-th Calabi flow, we generalize Ge-Hua's a theorem ( Theorem 3.4 in \cite{GB}) and give the following result, which show that the solution $r(t)$ to the combinatorial $p$-th Calabi flow (\ref{p-Cah}) in $H^{2}$ remains bounded away from 0 and bounded from above as well, as long as time remains bounded.
\begin{theorem}\label{bound-above}
\emph{Let $(M,\mathcal{T})$ be a triangulated compact surface with an edge weight $\Phi:E\rightarrow[0,\pi/2]$. Let $r(t)$ be the unique solution to the combinatorial $p$-th Calabi flow (\ref{p-Cah}) in hyperbolic background geometry on the maximal time interval $[0,T)$. Then all $r_i(t)$ are uniformly bounded above on $[0,T)$}.
\end{theorem}
\noindent\begin{proof}
The result is arrived at by contradiction. The proof is similar to the proof of Theorem 3.4 in \cite{GB}.
Assume that it is not true, then there exists at least one vertex $i\in V$, such that
\begin{equation}
\limsup_{t\rightarrow T}r_i(t)=+\infty.
\end{equation}
For this vertex $i$, using Lemma 3.5 in \cite{CL}, we can choose a sufficient large positive number $l$ such that $r_i>l$ and the inner angle $\theta_i$ is smaller than $\frac{\pi}{2d_i}$, where $d_i$ is the degree of the vertex $i$. Then we have $K_i>\pi$.\par
 Set $L=\max\{l,C,r_i(0)+1\}$, where $C$ is given in Proposition \ref{Area-e}. We claim that if $r_i(t)>L$, then
\begin{equation}\label{<0}
\frac{dr_i}{dt}<0
\end{equation}
 for any $t\in(0,T)$. To see this, first by the definition of the combinatorial curvature, i.e., $K_i=2\pi-\sum_{{ijk}\in F}\theta_i^{jk}$, we have $K_i<2\pi$.\par
  Combining above inequlities on curvature, we have $K_j-K_i<\pi$ for $j\thicksim i$. Note that if $a\leq b$, then $|a|^{p-2}a\leq|b|^{p-2}b$. Consequenly, by the $p$-th combinatorial Calabi flow (\ref{p-Cah}) and Proposition \ref{Area-e} for $\epsilon=\frac{1}{\pi^{p-2}}$, we have
\begin{align*}
\frac{1}{\sinh r_i}\frac{dr_i}{dt}&=\sum_{j\thicksim i}B_{ij}|K_j-K_i|^{p-2}(K_j-K_i)-A_iK_i\\&\leq\sum_{j\thicksim i}B_{ij}\pi^{p-1}-A_i\pi\\&=\pi^{p-1}\left(\sum_{j\thicksim i}B_{ij}-\frac{1}{\pi^{p-2}}A_i\right)\\&\leq0
\end{align*}
This proves the claim.

We may choose $t_0\in(0,T)$ such that $r_i(t_0)>c$ by (\ref{<0}). Suppose that $t_1\in[0,t_0]$ attains the maximum of $r_i(t)$ in $[0,t_0]$. Then by the definition of $L$, $t_1>0$. Thus
$$\frac{d r_i}{d t}(t_1)\geq 0,$$
which contradicts to (\ref{<0}). This completes the proof. $\hfill\Box$

Now we can prove Theorem \ref{p-ca-exis} as follows. \\
\noindent\begin{proof}
In the case of $E^{2}$, the proof is similar to the proof of Theorem 3.5 in \cite{Ge4}.
Let $d_{i}$ denote the degree at vertex $v_{i}$, which is the number of edges adjacent
to $v_{i}$. Set $d=\max\{d_1,\ldots,d_N\}$, then $(2-d)\pi < K_i < 2\pi$, hence
$$|K_{j}-K_{i}|<d\pi,\ \forall j\sim i.$$ By the estimation
of $B_{ij}$ in Proposition \ref{B}, all $|\Delta_{p}K_{i}|$ are uniformly bounded by a positive constant
$c=2\sqrt{3}\cdot N\cdot (d\pi)^{p-1}$, which depends only on the triangulation. Then we have
$$c_{0}e^{-ct}\leq r_{i}(t)\leq c_{0}e^{ct}$$
where $c_{0} = c(r(0))$, which implies that the combinatorial $p$-th Calabi flow (\ref{p-Ca}) has a solution
for all time $t\in [0, +\infty)$ for any $r(0)\in \mathbb{R}^{N}_{>0}.$
In the case of $H^{2}$, Theorem \ref{bound-above} show that as long as time remains bounded, then the solution $r_{i}(t)$ to the $p$-th combinatorial Calabi flow (\ref{p-Cah}) is bounded away from 0 and bounded from above. This implies that
the solution $r_{i}(t)$ to the $p$-th combinatorial Calabi flow (\ref{p-Cah}) exists for all time $t\geq0$.
 $\hfill\Box$
\end{proof}

\section{The long time convergence}
\hspace{14pt}  In this section we prove the convergence parts of both Theorem \ref{Ca-E} and Theorem \ref{Ca-H}. More precisely, we show the long time convergence of both the solutions to the combinatorial $p$-th Calabi flows (\ref{p-Ca}) in $E^{2}$ and (\ref{p-Cah}) in $H^{2}$. In fact, we have
\begin{theorem}\label{Ca-conv}
\emph{Suppose $\{r(t)|t\in[0, +\infty)\}$ is a long time solution to the combinatorial $p$-th Calabi flow \ref{p-Cah} in hyperbolic background geometry, then $\{r(t)|t\in[0, +\infty)\}$ converges if and only if there exists a zero curvature circle packing metric $r_{ze}$.
Suppose $\{r(t)|t\in[0, +\infty)\}$ is a long time solution to the combinatorial $p$-th Calabi flow \ref{p-Ca} in Euclidean background geometry with any initial value $r(0)$ satisfying $\prod_{i=1}^{N}r_{i}(0)=1$, then $\{r(t)|t\in[0, +\infty)\}$ converges if and only if there exists a constant curvature circle packing metric $r_{av}$}.
\end{theorem}
First we give the following compactness theorem which is key to prove Theorem \ref{Ca-conv} .
\begin{proposition}\label{conv-he}
\emph{Suppose $u(t)$ for $t\in[0, +\infty)$ is a solution to the combinatorial $p$-th Calabi flow (\ref{p-Cah}) (\ref{p-Ca} resp.) in the hyperbolic (Euclidean resp.) background geometry so that the set $\{u(t)\mid t\in[0, +\infty)\}$ lies in a compact region in $\mathbb{R}^{N}$. Then $u(t)$ converges a zero (constant resp.) curvature circle packing metric $u_{ze}$ ($u_{av}$ resp.)}.
\end{proposition}

\subsection{Some properties on the combinatorial $p$-Laplace operator}
In this subsection we show some properties with respect to the combinatorial $p$-Laplace operator $\Delta_p$ (\ref{p-Lap}). \par
For any smooth closed manifold $(M,g)$ with Riemannian metric $g$, considering the smooth $p$-Laplace operator
$\Delta_pf=\text{div}(|\nabla f|^{p-2}\nabla f),$ where $f$ is an arbitrary smooth function on $M$, and $\text{div}$ denotes the divergence. By the divergence theorem, we have
$\oint_{M}\Delta_pfd\mu=0.$
For a combinatorial
surface $(X, T, \Phi)$ with a circle packing metric $r$, similarly we have
\begin{lemma}\label{Div}
\emph{For any $f\in C(V)$, we have
$$\sum_{i}^{N}\Delta_pf_{i}=0.$$
Specially, for the curvature function $K_{i}$ we have}
$$\sum_{i=1}^{N}\Delta_pK_{i}=0.$$
\end{lemma}

By above lemma and the combinatorial $p$-th Calabi flow equation \ref{p-Ca-u} we generalize Ge's Proposition 3.4 in \cite{Ge4} and obtain the following important invariance property.
\begin{proposition} \label{inv}
\emph{As long as the solution $u(t)$ (or equivalently $r(t)$) to the combinatorial $p$-th Calabi flow \ref{p-Ca-u} (or equivalently \ref{p-Ca}) in $E^{2}$ exists, both
$\prod_{i=1}^{N} r_{i}(t)\equiv\prod_{i=1}^{N} r_{i}(0)$ and
$\sum_{i=1}^{N} u_{i}(t)\equiv\sum_{i=1}^{N} u_{i}(0)$ are constants}.
\end{proposition}
\begin{lemma}\label{iden}
\emph{For the $p$-Laplace operator $\Delta_p$ and For any $f\in C(V)$, we have the following identity } $$\sum_{i=1}^{N}f_i\Delta_pf_{i}=-\frac{1}{2}\sum_{i=1}^{N}\sum_{j\thicksim i}B_{ij}|f_j-f_i|^{p}.$$
\emph{Specially, for the curvature function $K_{i}$ we have}
$$\sum_{i=1}^{N}K_i\Delta_pK_{i}=-\frac{1}{2}\sum_{i=1}^{N}\sum_{j\thicksim i}B_{ij}|K_j-K_i|^{p}.$$
\end{lemma}
\noindent\begin{proof}
By the definition (\ref{p-Lap}) of $\Delta_p$ and the symmetry of $B_{ij}$, we have
\begin{align*}\sum_{i=1}^{N}f_i\Delta_pf_{i}&=\sum_{i=1}^{N}\sum_{j\thicksim i}B_{ij}|f_j-f_i|^{p-2}(f_j-f_i)f_i\\
&=\frac{1}{2}\sum_{i=1}^{N}\sum_{j\thicksim i}B_{ij}|f_j-f_i|^{p-2}(f_j-f_i)f_i\\
&+\frac{1}{2}\sum_{i=1}^{N}\sum_{j\thicksim i}B_{ij}|f_j-f_i|^{p-2}(f_j-f_i)(-f_j)\\
&=-\frac{1}{2}\sum_{i=1}^{N}\sum_{j\thicksim i}B_{ij}|f_j-f_i|^{p}.
\end{align*}
This completes the proof.  $\hfill\Box$
\end{proof}
For any $c > 0$, denote $\mathcal{P}_{c} =\{r = (r_{1}, . . . , r_{N})^{T} \in \mathbb{R}^{N}_{>0}
\mid \prod_{i=1}^{N} r_{i}=c\}$ and denote $\mathcal{U}_{a}=\{u = (u_{1}, . . . , u_{N})^{T}\mid\sum_{i=1}^{N} u_{i}=a\}$. Then  $\mathcal{P}_{c}=Exp(\mathcal{U}_{a})$, where $a=\text{ln} c$. By Proposition \ref{inv}, $\{r(t)\}\subset \mathcal{P}_{c}$ and $\{u(t)\}\subset \mathcal{U}_{a}$ along the combinatorial $p$-th
Calabi flow in $E^{2}$, where $c =\prod_{i=1}^{N} r_{i}(0)$ and $a=\sum_{i=1}^{N} u_{i}(0)$.\par

\subsection{Proof of Proposition \ref{conv-he}}
Now we give the proof of Proposition \ref{conv-he}.\\
\noindent\begin{proof}
In the case of $H^{2}$, suppose $u(t)$ for $t\in[0, +\infty)$ is a solution to the combinatorial $p$-th Calabi flow (\ref{p-Cah-u}).
By the assumption, we know $\{u(t)\mid t\in[0, +\infty)\}\subset M$, where $M$ is a compact subset in $\mathbb{R}^{N}$. Consider the combinatorial Ricci potential
\begin{equation*}
F(u)=\int_{u_{ze}}^u\sum_{i=1}^{N}K_idu_i,\ u\in\mathbb{R}_{<0}^{N},
\end{equation*}
where $u_{ze}=\text{ln} r_{ze}$. \par
Set $\varphi(t)=F(u(t))$, then we carry out the proof in two steps.\\
\noindent\textbf{Step 1.} The solution $u(t)$  to the combinatorial $p$-th Calabi flow (\ref{p-Cah-u}) sequently converges to a zero curvature circle packing metric $u_{ze}$. \\
By Proposition \ref{rigidity}, the gradient map of the combinatorial Ricci potential
\begin{align*}
\nabla F: \mathbb{R}_{<0}^{N}&\rightarrow \mathbb{R}^{N}\\
u&\mapsto K
\end{align*}
is injective, and $u=u_{ze}$ is the only critical point. Furthermore, the Hessian of $F$ is exactly $L=(L_{ij})$, which is positive definite by Proposition \ref{L2}. Therefore, $F(u)\geq F(u_{ze})=0$.\par

Note that $B_{ij}=B_{ji}$, and by Lemma \ref{iden} we have
\begin{align*}
K^{T}\Delta_pK=-\frac{1}{2}\sum_{i=1}^{N}\sum_{j\thicksim i}B_{ij}|K_j-K_i|^{p}\leq0.
\end{align*}
Moreover, $A=diag\{A_1,A_2,\ldots,A_{N}\}$ is positive definite. Using the $p$-th combinatorial Calabi flow equation (\ref{p-Cah-u}), we can compute the derivative of $\varphi(t)$ as follows
\begin{equation*}
\varphi'(t)=\sum_{i}K_iL_pK_i=K^TL_pK=K^T(\Delta_p-A)K\leq 0,
\end{equation*}
which implies that $\varphi(t)$ is decreasing and bounded, i.e.,
$0\leq F(u(t))=\varphi(t)\leq \varphi(0)=F(u(0))$. Hence $\varphi(t)$ converges, i.e., when $t\rightarrow +\infty$,
$\varphi(t)\rightarrow \varphi(+\infty)$.
By the mean value theorem, there exists $t_{n}\in(n, n+1)$ satisfying
$\varphi(n+1)-\varphi(n)=\varphi'(t_{n})$, thus $\varphi'(t_{n})\rightarrow 0$ when $n\rightarrow +\infty$.

By the assumption,  we have ${u(t_{n})}\subset M$,
which means that there exists a convergent subsequence $t_{n_{k}}\uparrow +\infty$ s.t.
$\lim_{k\rightarrow+\infty}u(t_{n_{k}})=u^{*}.$ Consider $K_{i}, B_{ij}$  as smooth functions of $u$, which means there exists $ K_{i}^{*}, B_{ij}^{*}$ s.t.
$\lim_{k\rightarrow+\infty}K_{i}(t_{n_{k}})=K_{i}^{*}$ and $\lim_{k\rightarrow+\infty}B_{ij}(t_{n_{k}})=B_{ij}^{*}.$
By equalities above $\varphi'(t)=K^T(\Delta_p-A)K$ and $$K^{T}\Delta_pK=-\frac{1}{2}\sum_{i=1}^{N}\sum_{j\thicksim i}B_{ij}|K_j-K_i|^{p}\leq0,$$ we have $$0\geq-\frac{1}{2}\sum_{i=1}^{N}\sum_{j\thicksim i}B_{ij}^{*}|K_j^{*}-K_i^{*}|^{p}=(K^{*})^TA^{*}K^{*}\geq0,$$
which implies $(K^{*})^TA^{*}K^{*}=0$. Note that $A=diag\{A_1,A_2,\ldots,A_{N}\}$ is always positive definite on the compact subset $M$, hence $A^{*}$ is positive definite, which implies that $K_{i}^{*}=0$, for $i=1, \cdots, N$, hence $u^{*}$ is exactly a zero curvature circle packing metric $u_{ze}$, i.e.,
$$\lim_{k\rightarrow+\infty}u(t_{n_{k}})=u_{ze}.$$ \\

\noindent\textbf{Step 2.}  $u(t)\rightarrow u_{ze}, t\rightarrow+\infty.$\\
We claim \begin{equation}\label{limit}
\lim_{t\rightarrow+\infty}\varphi(t)=0.
\end{equation}
Recall $\varphi(t)$ converges, i.e., when $t\rightarrow +\infty$,
$\varphi(t)\rightarrow \varphi(+\infty)$.
By $\lim_{k\rightarrow+\infty}u(t_{n_{k}})=u_{ze}$,
we have $$\varphi(t_{n_{k}})=F(u(t_{n_{k}}))\rightarrow F(u_{ze})=0, k\rightarrow +\infty,$$
therefore $\varphi(+\infty)=0$, i.e.,
\begin{equation*}
\lim_{t\rightarrow+\infty}\varphi(t)=0.
\end{equation*}
Now we prove $$\lim_{t\rightarrow+\infty}u(t)=u_{ze}$$ in fact by contradiction. If not, then there exists $\delta>0$ and $\xi_{n}\uparrow +\infty$ such that
$$|u(\xi_{n})-u_{av}|>\delta,$$
 which means $\{u(\xi_{n})\}\subset \mathbb{R}^{N}\setminus B(u_{av}, \delta)$.  Recall $\{u(t)\mid t\in [0, +\infty)\}\subset M$ where $M$ is a compact subset of $\mathbb{R}^{N}$. Hence $\{u(\xi_{n})\}\subset M\cap (\mathbb{R}^{N}\setminus B(u_{av}, \delta))$. It is obvious that continuous differentiable function $F(u)$ has a positive lower bound in the compact set $M\cap (\mathbb{R}^{N}\setminus B(u_{av}, \delta))$, i.e.,
$$F(u)\geq c>0, \forall u\in M\cap (\mathbb{R}^{N}\setminus B(u_{av}, \delta)).$$
Hence we have $$\varphi(\xi_{n})= F(u(\xi_{n}))\geq c>0,$$
which causes a contradiction considering $\varphi(\xi_{n})\rightarrow0, n\rightarrow +\infty$ by (\ref{limit}).
Therefore $$\lim_{t\rightarrow+\infty}u(t)=u_{ze}.$$

In the case of $E^{2}$, consider the combinatorial Ricci potential
\begin{equation*}
F(u)=\int_{u_{av}}^u\sum_{i=1}^{N}(K_i-k_{av})du_i,\ u\in\mathbb{R}^{N},
\end{equation*}
where $u_{av}=\text{ln} r_{av}$. \par
By Proposition \ref{rigidity}, the gradient map
\begin{align*}
\nabla F\mid_{\mathcal{U}_{a}}:\mathcal{U}_{a}&\rightarrow \mathbb{R}^{N}\\
u&\mapsto K-K_{av}
\end{align*}
is injective, where $a=\sum_{i=1}^{N}u_{i}(0)$, and $u=u_{av}$ is the only critical point. Furthermore, the Hessian of $F$ is exactly $L=(L_{ij})$, which is semi-positive definite by Proposition \ref{L2}. Therefore, $F(u)\geq F(u_{av})=0$.\par
 Set $\varphi(t)=F(u(t))$, then we
compute the derivative of $\varphi(t)=F(u(t))$ as follows
\begin{align*}
\varphi'(t)&=\sum_{i=1}^{N}(K_i-k_{av})\Delta_pK_{i}\\
&=\sum_{i=1}^{N}K_i\Delta_pK_{i} \\
&=-\frac{1}{2}\sum_{i=1}^{N}\sum_{j\thicksim i}B_{ij}|K_j-K_i|^{p} \\&\leq0,
\end{align*}
where the second equality is due to Lemma \ref{Div}, the third equality is due to Lemma \ref{iden} and  the last inequality is due to Proposition \ref{B}. \par
This implies that $\varphi(t)$ is decreasing and bounded, i.e., $0\leq \varphi(t)\leq \varphi(0)=F(u(0))$. Then by the same arguments above in the case of $H^{2}$, we can show that $u(t)$ converges a constant curvature circle packing metric $u_{av}$.
 $\hfill\Box$
\end{proof}

\subsection{Proof of Theorem \ref{Ca-conv}}

Finally, we give the proof of Theorem \ref{Ca-conv}.\\
\noindent\begin{proof}
First we show that there exists a zero curvature circle packing metric $r_{ze}$ if the long time solution $\{r(t)|t\in[0, +\infty)\}$ or equivalently $\{u(t)|t\in[0, +\infty)\}$ to the combinatorial $p$-th Calabi flow \ref{p-Cah-u} in $H^{2}$ converges, i.e.,
$$u^{*}=\lim_{t\rightarrow+\infty}u(t), u_{i}^{*}=\lim_{t\rightarrow+\infty}u_{i}(t)$$ exists.
Hence by the mean value theorem, there exists $t_{n}\in(n, n+1)$ satisfying
$$u_{i}(n+1)-u_{_{i}}(n)=u'_{i}(t_{n}).$$
It is obvious that the limit of the left of above equation is $0$ when $n\rightarrow +\infty$, so by the combinatorial $p$-th Calabi flow equation (\ref{p-Cah-u}) we have
$$L_pK_i\mid_{u=u^{*}}=(\Delta_p-A_i)K_i\mid_{u=u^{*}}=\lim_{n\rightarrow+\infty}(\Delta_p-A_i)K_i\mid_{u=u(t_{n})}=\lim_{n\rightarrow+\infty}u'_{i}(t_{n})=0.$$
Therefore $\sum_{i}K_iL_pK_i\mid_{u=u^{*}}=K^T(\Delta_p-A)K\mid_{u=u^{*}}=0$, by the same arguments in the proof above, we have
$$\lim_{n\rightarrow+\infty}K_{i}(t_{n})=K_{i}(u^{*})=0, i=1, \cdots, N.$$
This means $u^{*}$ is actually a zero curvature circle packing metric $u_{ze}$.
Using the same argument we can show that there exists a constant curvature circle packing metric $r_{av}$ if the long time solution $\{r(t)|t\in[0, +\infty)\}$ or equivalently $\{u(t)|t\in[0, +\infty)\}$ to the combinatorial $p$-th Calabi flow (\ref{p-Ca-u}) in $E^{2}$ converges.\par
 As for the  " if " part of the proof, we divide it into cases:\par
In the case of $H^{2}$, suppose $\{r(t)|t\in[0, +\infty)\}$ or equivalently $\{u(t)|t\in[0, +\infty)\}$ is a solution to the combinatorial $p$-th Calabi flow (\ref{p-Cah-u}) in $H^{2}$. If there exists a zero curvature circle packing metric $r_{ze}$ or equivalently $u_{ze}$, consider the the combinatorial Ricci potential $F$ (\ref{RP1}), and set $\varphi(t)=F(u(t))$. We know that $\varphi(t)$ is decreasing and bounded, i.e., $0\leq\varphi(t)\leq \varphi(0)=F(u(0))$ by the proof above.  By Theorem \ref{proper}, the combinatorial Ricci potential $F$ (\ref{RP1}) is proper,  thus $F^{-1}[0, F(u(0))]$ is a compact subset ,which implies $\{u(t)\mid t\in[0, +\infty)\}$ lies in a compact region $F^{-1}[0, F(u(0))]$ in $\mathbb{R}^{N}$. By Proposition \ref{conv-he}, $\{u(t)\mid t\in[0, +\infty)\}$ converges the zero curvature circle packing metric $u_{ze}$. \par

In the case of $E^{2}$,  if there exists a constant curvature circle packing metric $r_{av}$ or equivalently $u_{av}$, and suppose $u(t)$ for $t\in[0, +\infty)$ is a solution to the combinatorial $p$-th Calabi flow (\ref{p-Ca-u}) in $E^{2}$ with any initial value $\sum_{i=1}^{N}u_{i}(0)=0$ (equivalently $\prod_{i=1}^{N}r_{i}(0)=1$). Note by Proposition \ref{inv}, $\{u(t)\mid t\in [0, +\infty)\}\subset  \mathcal{U}_{a}$ along the combinatorial $p$-th
Calabi flow (\ref{p-Ca-u}), where $a=\sum_{i=1}^{N}u_{i}(0)=0$ in this case. Then consider the combinatorial Ricci potential
\begin{equation*}
F(u)=\int_{u_{ze}}^u\sum_{i=1}^{N}(K_i-k_{av})du_i,\ u\in\mathbb{R}^{N},
\end{equation*} and set $\varphi(t)=F(u(t))$. By the proof above we know that $\varphi(t)$ is decreasing and bounded, i.e., $0\leq\varphi(t)\leq \varphi(0)=F(u(0))$ . Hence $\{u(t)\mid t\in [0, +\infty)\}\subset F\mid_{\mathcal{U}_{0}}^{-1}[0, \varphi(0)]$. By Proposition \ref{proper} $F\mid_{\mathcal{U}_{0}}$ is proper, thus $F\mid_{\mathcal{U}_{0}}^{-1}[0, \varphi(0)]$ is a compact subset of $\mathcal{U}_{0}$, which implies $\{u(t)\mid t\in[0, +\infty)\}$ lies in a compact region $F^{-1}[0, F(u(0))]$ in $\mathbb{R}^{N}$. By Proposition \ref{conv-he}, $\{u(t)\mid t\in[0, +\infty)\}$ converges the constant curvature circle packing metric $u_{av}$.
 $\hfill\Box$
\end{proof}

\section{Further discussions}
\hspace{14pt} Theorem \ref{Ca-E} and Theorem \ref{Ca-H} generalize the work of H. Ge, Ge-Xu and Ge-Hua from the combinatorial Calabi flow to the combinatorial $p$-th Calabi flow. Moreover, our results show that the combinatorial $p$-th Calabi flow does not have the exponential convergence in general $p>1$. There are still several interesting questions related to our results.

\subsection{Normalized combinatorial $p$-th Calabi flow}
\hspace{14pt} The \emph{normalized combinatorial $p$-th Calabi flow in Euclidean background geometry} is defined as follows.
  \begin{equation}\label{p nc-flow}
  u'_{i}(t)=-\nabla_{u_{i}}E_{p}=p\sum_{j=1}^{N}L_{ij}|k_{av}-K_{j}|^{p-2}(k_{av}-K_{j}),
\end{equation}
where $E_{p}(u)$ is the $p$-th Calabi energy, i.e., $E_{p}(u)=\sum_{i=1}^{N}|k_{av}-K_{i}|^{p}$.\par
It is obvious that (\ref{p nc-flow}) is precisely the negative gradient flow of the $p$-th Calabi energy $E_{p}$.
We can similarly show the long time existence of (\ref{p nc-flow}).  If there exists the compactness of the closure of the solution subset $\{u(t)\mid t\in[0, +\infty)\}$, then we can show the long time convergence for the curvature flow  $\{K(u(t))\mid t\in[0, +\infty)\}$ when $t\rightarrow +\infty$. To see this, by the compactness of the closure of the solution subset $\{u(t)\mid t\in[0, +\infty)\}$ and the similar argument of the proof of Proposition \ref{conv-he}, we can show that the solution $u(t)$ to (\ref{p nc-flow}) sequently converges a constant curvature circle packing metric $u_{av}$. That is to say, there exists a sequence $\xi_{n}\rightarrow +\infty$ such that

$$\lim_{n\rightarrow+\infty}u(\xi_{n})=u_{av}.$$

 By (\ref{p nc-flow}) we have
\begin{align*}
\frac{\partial E_{p}(u(t))}{\partial t}
=\sum_{i=1}^{N}\nabla_{u_{i}}E_{p}u'_{i}
=-\sum_{i=1}^{N}(u'_{i})^{2}
\leq 0,
\end{align*}
which implies that $E_{p}(u(t))$ is decreasing and bounded, i.e., $0\leq E_{p}(u(t))\leq E(u(0))$.  Hence
$$\lim_{t\rightarrow+\infty}E_{p}(u(t))=A$$ exists.
In fact $A=0$. Note that
$$\lim_{n\rightarrow+\infty}u(\xi_{n})=u_{av},$$
thus
$$A=\lim_{t\rightarrow+\infty}E_{p}(u(\xi_{n}))=E_{p}(u_{av})=0.$$
Hence we have $$\lim_{t\rightarrow+\infty}E_{p}(u(t))=0.$$
Further we have
$$ 0\leq|k_{av}-K_{i}(u(t))|^{p}\leq E_{p}(u(t)), \forall i=1,\cdots, N$$
thus $$\lim_{t\rightarrow+\infty}K_{i}(u(t))=k_{av}, \forall i=1,\cdots, N.$$
By the above argument, the key question is to show the compactness of the closure of the solution subset $\{u(t)\mid t\in[0, +\infty)\}$ to (\ref{p nc-flow}). Or can one directly show the long time convergence of the solution $\{u(t)\mid t\in[0, +\infty)\}$ to (\ref{p nc-flow}) when $t\rightarrow +\infty$ using different approaches?

\subsection{Another combinatorial $p$-th Calabi type flow}
\hspace{14pt} Let $G = (V,E)$ be a finite graph, where $V$ denotes the vertex set and $E$ denotes the edge
set. Fix a vertex measure $\mu : V\rightarrow(0,+\infty)$ s.t., $\mu_{i}=\mu(i)\equiv1, i=1, \cdots, N$ and an edge measure $\omega : E \rightarrow(0,+\infty)$ on $G$.
The edge measure ¦Ø is assumed to be symmetric, that is, $\omega_{ij} = \omega_{ji}$ for each edge $i\sim j$. Then we study another form of the combinatorial $p$-th curvature flow. For any given combinatorial $p$-th Laplace operator, namely fix
$$ \Delta_{p}f_{i}=\frac{1}{\mu_{i}}\sum_{j\sim i}\omega_{ij}|f_{j}-f_{i}|^{p-2}(f_{j}-f_{i})=\sum_{j\sim i}\omega_{ij}|f_{j}-f_{i}|^{p-2}(f_{j}-f_{i}),$$
where $\mu_{i}=1$ and $\omega_{ij}\geq 0$ are independent on time $t$.\par
Consider the combinatorial $p$-th Calabi type flow
\begin{equation}\label{p-flow}
u'_{i}=\frac{du_{i}}{dt}=\Delta_{p}K_{i}.
\end{equation}
For these equations, we can similarly prove the long time existence.
Then consider the $p$-th energy $$E(u)=\int_{E}|\nabla K|^{p}d\omega=\sum_{j\sim i}\omega_{ij}|K_{j}-K_{i}|^{p},$$
and we have $$\frac{\partial E}{\partial K_{j}}=-p\Delta_{p}K_{j}.$$
Hence
\begin{align*}
\frac{dE(u(t))}{dt}=\sum_{i=1}^{N}\frac{\partial E}{\partial u_{i}}u'_{i}
=\sum_{i=1}^{N}(\sum_{j=1}^{N}\frac{\partial E}{\partial K_{j}}\frac{\partial K_{j}}{\partial u_{i}})u'_{i}
=-p\sum_{i=1}^{N}\sum_{j=1}^{N}u'_{i}L_{ij}u'_{j}\leq 0,
\end{align*}
where the last inequality holds because the matrix $L=(L_{ij})\geq 0$. \par
Therefore $E(u(t)$ is decreasing and bounded.
If there exists the compactness of the closure of the solution subset $\{u(t)\mid t\in[0, +\infty)\}$ to (\ref{p-flow}), then we can show the long time convergence for the curvature flow  $\{K(u(t))\mid t\in[0, +\infty)\}$ when $t\rightarrow +\infty$ using the same argument in above subsection. To see this, we first similarly show
$$\lim_{t\rightarrow+\infty}E(u(t))=0.$$
Next we we prove $$\lim_{t\rightarrow+\infty}K_{i}(t)=K_{av}.$$
For this we need a lemma due to H. Ge \cite{GX}
\begin{lemma}\label{P-ineqality}
Let $G = (V,E)$ be a finite graph. For all functions $\varphi\in C(V )$ with $\bar{\varphi}= 0$,
there exists some positive constant $C_{p,G}$ depending on $G$ and $p$ such that
\begin{equation}
\int_{V}|\varphi|^{p}d\mu\leq C_{p,G}\int_{E}|\nabla \varphi|^{p}d\omega.
\end{equation}
\end{lemma}
Let $\varphi=K_{av}-K$, then $\bar{\varphi}=0$ because
\begin{align*}
\bar{\varphi}=\sum_{i=1}^{N}(k_{av}-K_{i})\mu_{i}/Vol(G)=\sum_{i=1}^{N}(k_{av}-K_{i})/Vol(G)=2\pi\chi(X)-\sum_{i=1}^{N}K_{i}=0,
\end{align*}
where we use the combinatorial Gauss-Bonnet formula in the case of $E^{2}$.\par
By the lemma above, for $i=1, \cdots, N $ we have
\begin{align*}
|k_{av}-K_{i}|^{p}\leq\sum_{i=1}^{N}|k_{av}-K_{i}|^{p}
=\int_{V}|\varphi|^{p}d\mu
\leq C_{p,G}\int_{E}|\nabla \varphi|^{p}d\omega
=E(u(t)),
\end{align*}
which implies  $$\lim_{t\rightarrow+\infty}K_{i}(t)=k_{av}, \forall i,\cdots, N.$$
Therefore the key question is also to show the compactness of the closure of the solution subset $\{u(t)\mid t\in[0, +\infty)\}$ to (\ref{p-flow}). Or can one directly show the long time convergence of the solution $\{(u(t)\mid t\in[0, +\infty)\}$ to (\ref{p-flow}) when $t\rightarrow +\infty$ using different approaches?

\end{proof}

(Aijin Lin) Department of Mathematics, National University of Defense Technology, Changsha 410073, P. R. China\\
\emph{E-mail address}: \verb"linaijin@nudt.edu.cn"\\
\\
(Xiaoxiao Zhang) Institute of Mathematics, Beijing Jiaotong University, Beijing 100044, P. R. China\\
\emph{E-mail address}: \verb"xiaoxiaozhang0408@bjtu.edu.cn"\\
\end{document}